\documentclass[11pt]{article}
\newcommand{\documentdate}{16 VIII 2024}
\usepackage{a4wide,latexsym,graphicx,amsmath}
\title{Refining asymptotic complexity bounds for nonconvex optimization methods,
  including why steepest descent is $o(\epsilon^{-2})$ rather than $\calO(\epsilon^{-2})$}

\author{
S. Gratton%
  \thanks{Universit\'{e} de Toulouse, INP, IRIT, Toulouse, France. Email:
     serge.gratton@enseeiht.fr. Work partially supported by 3IA Artificial and
     Natural Intelligence Toulouse Institute (ANITI), French "Investing for the Future
     - PIA3" program under the Grant agreement ANR-19-PI3A-0004"}
~and C.-K. Sim%
  \thanks{University of Portsmouth, Hampshire, United Kingdom. Email: chee-khian.sim@port.ac.uk }
~and Ph. L. Toint%
  \thanks{NAXYS, University of Namur, Namur, Belgium. Email: philippe.toint@unamur.be}
}

\newcommand{\beqn}[1]{\begin{equation}\label{#1}}
\newcommand{\eeqn}{\end{equation}}
\newcommand{\req}[1]{(\ref{#1})}
\newcommand{\ms}{\;\;\;\;}
\newcommand{\tim}[1]{\;\; \mbox{#1} \;\;}

\newtheorem{theorem}{Theorem}[section]
\newcommand{\numsection}[1]{\section{#1}\setcounter{equation}{0}}

\newcommand{\iiz}[1]{\{ 0, \ldots, #1 \}}

\renewcommand{\Re}{\hbox{I\hskip -2pt R}}

\newcommand{\calS}{{\cal S}} 
\newcommand{\calO}{{\cal O}} 
\renewcommand{\Re}{\hbox{I\hskip -2pt R}}

\newcommand{\sfrac}[2]{{\scriptstyle \frac{#1}{#2}}}
\newcommand{\half}{\sfrac{1}{2}}
\newcommand{\comment}[1]{}

\renewcommand{\thefootnote}{(\arabic{footnote})}

\newcommand{\bpr}{{\bf Proof.} \hspace{1.5mm}}
\newcommand{\epr}{\hfill $\Box$ \vspace*{1em}}
\newcommand{\proof}[1]{
\renewcommand{\thefootnote}{\fnsymbol{footnote}}
\begin{list}{}{
\setlength{\topsep}{0.0pt}
\setlength{\partopsep}{0.0pt}
\setlength{\leftmargin}{0.025\textwidth}
\setlength{\rightmargin}{0.5\leftmargin}
\setlength{\labelwidth}{0.5\leftmargin}
\setlength{\labelsep}{0.25\leftmargin}}
\item \bpr #1 \epr \noindent
\end{list}}
\newcommand{\lthm}[2]{\vspace{\baselineskip} 
\noindent\framebox[\textwidth]{\parbox{0.95\textwidth}{
\begin{theorem} \label{#1} \rm #2 \end{theorem} } } \vspace{\baselineskip} }

\date{\documentdate}

\begin{document}

\maketitle

\begin{abstract}
We revisit the standard ``telescoping sum'' argument ubiquitous in the
final steps of analyzing evaluation complexity of algorithms for
smooth nonconvex optimization, and obtain a refined formulation of the
resulting bound as a function of the requested accuracy $\epsilon$. While bounds
obtained using the standard argument typically are of the form
$\calO(\epsilon^{-\alpha})$ for some positive $\alpha$, the refined
results are of the form $o(\epsilon^{-\alpha})$. We then explore to
which known algorithms our refined bounds are applicable and finally
describe an example showing how close the standard and refined bounds
can be.
\end{abstract}


\numsection{Introduction}

The numerical solution of nonlinear optimization problems often hinges
on descent algorithms, that is on algorithms in which a function (the
objective function, the residual, the merit function, etc.) is
monotonically decreasing in the course of the iterations. The analysis
of their iteration (and evaluation) complexity is then typically conducted
using a ``telescoping sum'' argument in which a lower bound of the
iteration-wise function decrease is summed in a ``telescoping sum''
over all iterations. Combining the resulting lower bounds with an upper bound on
the total decrease then yields an upper bound on the number of iterations
where the iteration-wise decrease is significant, in turn producing
the desired upper bound on the algorithm's worst-case behaviour.

Inspired by an unpublished note \cite{Sim21} of the second author on the steepest
descent method, the present paper revisits this telescoping argument,
which in turn results in an refined complexity bounds for a large
number of known optimization algorithms.

We first describe the refined argument in Section~\ref{arg-s}
and then investigate to which algorithms the new result is applicable
in Section~\ref{algos-s}.  We conclude by presenting, in
Section~\ref{ex-s}, an example showing that the new complexity bounds
may be very close to the standard ones.

\noindent
{\bf Notation.} Given two functions $a(\epsilon)$ and $b(\epsilon)$
with $b(\epsilon) > 0$ depending on a common parameter $\epsilon$
tending to zero, we say that $a(\epsilon) =
\calO\big(b(\epsilon)\big)$ if and only if there exists a constant
$\kappa < +\infty$ such that $\limsup_{\epsilon\rightarrow 0}
\big(a(\epsilon)/b(\epsilon)\big) \leq \kappa$.  We say that
$a(\epsilon) = o\big(b(\epsilon)\big)$ when this limit holds with
$\kappa = 0$. If $\calS$ is a set, $|\calS|$ denotes its
cardinality. Finally, $\lambda_{\min}(A)$ denotes the left-most
eigenvalue of the symmetric matrix $A$.

\numsection{A simple result about sequences}\label{arg-s}

In order to discuss our result, we need to consider the situation
where a specific optimization algorithm is applied to minimize a
smooth possibly nonconvex function
$f$ starting from $x_0$ and producing a sequence of iterates $\{x_k\}$, a sequence of
decreasing function values $\{f_k\}$ at these iterates and a sequence
of associated optimality measures $\{\omega_k\}$. We also need to
consider the set of ``successful iterations'' $\calS = \{k\geq 0 \mid x_{k+1} \neq x_k\}$.

Our results are asymptotic in the sense that we consider these sequences to be
infinite and examine how
\beqn{keps-def}
k(\epsilon) = \min\{ k \geq 0 \mid \omega_k \leq \epsilon \}
\eeqn
depends on $\epsilon$ when more and more accuracy
is requested, that is when $\epsilon$ tends to zero.
However, since the generation of these sequences will vary across the examples we will consider, we first
state our result in a slightly more abstract form.

\lthm{result}{
  Let $\{x_k\}$ be a sequence of iterates, $\{f_k=f(x_k)\}$ be a monotonically decreasing sequence bounded
  below, $\{\omega_k=\omega(x_k)\}$ be a non-negative sequence of
  optimality measures and let $\calS = \{k\geq 0 \mid x_{k+1} \neq
  x_k\}$ of successful iterations. Suppose also that
  \vspace*{-2mm}
  \beqn{kstop}
  \calS \cap \{k\geq 0 \mid \omega_k=0\} = \emptyset
  \eeqn
  and that
  \vspace*{-2mm}
  \beqn{succ}
  f_k-f_{k+1} \geq \kappa_d \, \omega_k^\beta
  \tim{ for} k \in \calS,
  \eeqn
  where $\kappa_d\in (0,1]$ and $\beta>0$ are constants. Suppose also there exist constants
  $\kappa_a\geq 1$ and $\kappa_b,\kappa_c \geq 0$ such that
  \beqn{growth}
  k \leq \kappa_a |\calS_k| + \kappa_b |\log(\omega_k)| + \kappa_c
  \tim{ whenever } \omega_k >0,
  \eeqn
  where $\calS_k=\calS \cap \iiz{k}$.
  Then,
  \vspace*{-5mm}
  \beqn{limomega}
  \lim_{k \rightarrow \infty} \omega_k = 0.
  \eeqn
  Moreover,
  either $k(\epsilon)$ is constant for all $\epsilon$ sufficiently small, or
  \beqn{bound2a}
  k(\epsilon)
  \leq \kappa_a\max\left[1,\frac{2(f_{\ell(k(\epsilon)-1)}-f_{k(\epsilon)})}{\kappa_d\epsilon^\beta}\right]
  +\kappa_b|\log(\epsilon)|+\kappa_a+\kappa_c,
  \eeqn
  where $\ell(k)$ is the largest index smaller or
  equal to the median of the indexes in $\calS_k$ and
  \beqn{limdiff}
  \lim_{\epsilon \rightarrow 0} \big(f_{\ell(k(\epsilon)-1)}-f_{k(\epsilon)-1}\big) = 0.
  \eeqn
  In all cases,
  \vspace*{-3mm}
  \beqn{bound2}
  k(\epsilon) = o\big(\epsilon^{-\beta}\big).
  \eeqn
}

  \proof{
    Suppose first that there exists a first $k_*\geq 0$ such that
    $\omega_{k_*} = 0$.  Then $k_*\not\in\calS$.  Thus $x_{k+1}=
    x_k$ and $\omega_{k+1}=0$, so that, by induction, $\omega_k=0$ for all $k\geq k_*$, implying
    \req{limomega}.  Moreover, $k(\epsilon) = k_*$ for all $\epsilon < \omega_{k_*-1}$.

    Suppose now that $\omega_k > 0$ for all $k\geq 0$.  If $|\calS|$
    is finite, then $\omega_k=\omega_{\min}$ for some
    $\omega_{\min}>0$ and
    all $k$ sufficiently large.  As a consequence and given that
    $|\calS_k|\leq|\calS|<+\infty$, the right-hand side
    of \req{growth} is bounded. But this is impossible since the
    left-hand side tends to infinity.  Hence $|\calS|$ is infinite.
    
    Consider now an arbitrary $k$ for which $|\calS_k|\geq 2$. Then $\ell(k)$ is well-defined
    and tends to infinity with $k$. We also have, using \req{succ} and
    the definition of $\ell(k)$, that
  \beqn{e1}
  f_{\ell(k)}-f_{k+1}
  = \sum_{j=\ell(k)}^k\big(f_j-f_{j+1}\big)
  = \sum_{j=\ell(k),j\in\calS_k}^k\big(f_j-f_{j+1}\big)
  \geq \half |\calS_k| \kappa_d \min_{j \in \calS_k}\omega_j^\beta.
  \eeqn
  Moreover, since $\{f_k\}$ is monotonically decreasing and bounded
  below, it is convergent, and hence
  \beqn{e2}
  \lim_{k\rightarrow \infty} \big(f_{\ell(k)}-f_{k+1}\big) = 0.
  \eeqn
  Thus the left-hand side of \req{succ} tends to zero with $k$, implying that
  $\lim_{k\rightarrow \infty,k\in\calS}\omega_k=0$. The
  definition of $\calS\supset\calS_k$ then ensures \req{limomega}.
  As a consequence $k(\epsilon)$ is well-defined for all $\epsilon$
  sufficiently small. By definition of $k(\epsilon)$ we also know that
  $\omega_k > \epsilon$ for all $k \leq k(\epsilon)-1$. Combining
  this inequality with \req{e1}, we obtain that
  \beqn{bound}
  |\calS_{k(\epsilon)-1}|\leq \frac{2(f_{\ell(k(\epsilon)-1)}-f_{k(\epsilon)})}{\kappa_d\epsilon^{\beta}}.
  \eeqn
Observe now that $k(\epsilon)$ is non-decreasing when $\epsilon$
tends to zero. Given its integer nature, either $k(\epsilon)$ tends to
infinity or is constant for all sufficiently small $\epsilon$.  In the
former case rewriting \req{e2} for $k=k(\epsilon)-1$ gives
\req{limdiff} and, because of \req{growth} taken at $k(\epsilon)$,
\beqn{nearfinal}
k(\epsilon)
\leq \kappa_a|\calS_{k(\epsilon)}|+\kappa_b \log(\omega_{k(\epsilon)}) + \kappa_c
\leq \kappa_a\big(|\calS_{k(\epsilon)-1}|+1) + \kappa_b|\log(\epsilon)|
+ \kappa_c,
\eeqn
which, given \req{bound}, yields \req{bound2a}.

We now prove \req{bound2}.
If $k(\epsilon)$ tends to infinity when $\epsilon$ tends to zero,
\req{bound2} is obtained by substituting \req{limdiff} in
\req{bound2a} and using the fact that $|\log(\epsilon)|= o(\epsilon^{-\beta})$.
If $k(\epsilon)$ remains constant, then \req{bound2}  immediately
follows from the fact that $\epsilon^{-\beta}$ tends to infinity when
$\epsilon$ goes to zero. 
} 

\noindent
Our assumption \req{kstop} simply says that, if, luckily, an exact
critical point of the desired order is found after finitely many
iterations, than the algorithm does not move away.  Note that we have
chosen $\ell(k)$ above to approximate the index separating $\calS_k$
in two parts of same cardinality, but other fixed proportions may of
course be used, at the price of modifying the constants in
\req{bound2a} and \req{bound}. Observe also that we could have
replaced $|\log(\omega_k)|$ in \req{growth} by any positive sequence
$\{h(x_k)\}$ such that $h(x_{k(\epsilon)}) =
o\big(\epsilon^{-\beta}\big)$.

\numsection{Application to existing algorithms and associated\\
  complexity bounds}\label{algos-s}

We now investigate the consequences of using this simple result
in the context where the sequence $\{x_k\}$ is the sequence of
iterates generated by specific nonlinear optimization algorithms
applied to sufficiently smooth functions that are bounded below. This
section only partially explores the resulting refined complexity
bounds, focusing on the algorithms described in the comprehensive book
\cite{CartGoulToin22}, but the authors are of course aware that the
discussion is incomplete.

\subsection{Unconstrained optimization}
\subsubsection{Steepest descent and other linesearch methods}

We start by considering complexity results for linesearch methods for
finding first-order critical points, such as those covering steepest
descent with Armijo, Goldstein \cite[Th.~2.2.2]{CartGoulToin22},
exact linesearch \cite[Th.~2.2.4]{CartGoulToin22} or with
Nesterov stepsize (\cite{Nest04} and
\cite[Equation~(2.2.5)]{CartGoulToin22}).  The proof of these results 
directly involves the ``telescoping sum'' argument, which we now cast
in the context of the previous section by selecting
\[
\{f_k\} =\{f(x_k)\},
\ms
\omega_k = \|\nabla_x^1f(x_k)\|
\ms
\beta = 2,
\ms
\calS_k = \iiz{k}
\]
and $\kappa_d$  is an algorithm-specific constant proportional to the
inverse of the gradient Lispchitz constant. Note that a standard
linesearch ensures that $\{f_k\}$ is decreasing in that \req{succ}
holds at all iterations. Morever the identity $\calS_k = \iiz{k}$
gives that $\kappa_a=1$ and $\kappa_b=\kappa_c=0$.

As a consequence, Theorem~\ref{result} implies that \emph{the
worst-case complexity of all these first-order algorithms (as a
function of the accuracy parameter $\epsilon$) is
$o\big(\epsilon^{-2}\big)$ rather than
$\calO\big(\epsilon^{-2}\big)$} as stated in the quoted
theorems. An illustration for steepest descent is discussed in
Section~\ref{ex-s}.

Interestingly, our technique does not require the complete sequence of
function values to satisfy \req{succ}, but it is
enough that these conditions hold, as is the case in the non-monotone ``gradient-related''
linesearch method discussed in \cite{CartSampToin15}, for a
subsequence of values at ``reference iterations'' which is used in the
telescoping sum argument.  Classical gradient-related linesearch
methods \cite{ShulSchnByrd85} are obtained by choosing the memory parameter
in this latter method to enforce monotonicity.

\subsubsection{Trust-region methods}

We may now turn to standard trust-region methods, whose complexity was
first considered in \cite{GratSartToin08} and is discussed in 
\cite[Th.~2.3.7 and 3.2.1]{CartGoulToin22} (for convergence of first- and
second-order methods converging to first-order critical points) and
\cite[Th.~3.2.6]{CartGoulToin22} for convergence to second-order
ones. Again the quoted proofs use a ``telescoping sum'' argument
where $\kappa_d$ an algorithm-specific constant proportional to the
inverse of the gradient Lispchitz constant
\[
\{f_k\} =\{f(x_k)\},
\ms
\omega_k = \|\nabla_x^1f(x_k)\| \tim{or}
\omega_k = \max\left[\|\nabla_x^1f(x_k)\|,\max(0,\lambda_{\min}(\nabla_x^2f(x_k))) \right],
\]
but we now choose $\calS_k$ to be the index set of the
``successful iterations'', that is iterations where $x_{k+1}$ differs
from $x_k$ and ensuring \req{succ}. The parameter
$\beta$ now depends on the purpose of the algorithm (finding first- or
higher-order critical points) and the degree of the objective's
derivatives used by the algorithm.  For standard trust-region methods
that seek first-order critical points, the parameter $\beta$ is
typically equal to two, while it is equal to three if second-order
ones are sought. Verifying \req{growth} is a little more complicated.
\cite[Lemma~2.3.1]{CartGoulToin22} shows that this inequality holds
with $\omega_k$ replaced by a lower bound on the trust-region radius.
Fortunately, \cite[Lem.~2.3.4 and 3.2.5]{CartGoulToin22} then state
that this lower bound is itself bounded below by $\omega_k$ or
$\epsilon$ (for $k=k(\epsilon)$), hence providing the desired inequality.

We may thus again apply our results to revisit all these
proofs. For the search of first-order critical points, this gives
\emph{$o\big(\epsilon^{-2}\big)$ rather than
$\calO\big(\epsilon^{-2}\big)$ complexity bounds} as a function of
$\epsilon$. The bounds for finding second-order points are similarly
refined to \emph{$o\big(\epsilon^{-3}\big)$ rather than
$\calO\big(\epsilon^{-3}\big)$}. In the
same vein, we may even consider trust-region methods for delivering
critical points of order higher than two
\cite[Th.~12.2.5]{CartGoulToin22} and obtain
$o\big(\epsilon^{-(q+1)}\big)$ worst-case complexity to
compute $q$-th order critical points. Finally, the global rates of
convergence of {\sf TR$q$IDA} and {\sf TR$q$EDA} trust-region variants for noisy
problems may also be refined in the same way (see \cite[Th.~13.1.8,
  13.3.4]{CartGoulToin22} together with
\cite[Lem.~13.1.1 and 13.1.4]{CartGoulToin22}).

But we may also consider more
elaborate trust-region-like algorithms, such as TR$\epsilon$
\cite{CurtRobiRoyeWrig21} (whose complexity proofs can be found in
\cite[Th.~3.4.5 and 3.4.6]{CartGoulToin22}), TRACE
\cite{CurtRobiSama18} (see \cite[Th.~3.4.11 and
  3.4.12]{CartGoulToin22} for proofs) or the Birgin-Martinez proposal
\cite{BirgMart17}. These methods achieve a complexity bound using $\beta = 3/2$ when
first-order points are sought.  Note that a specific result
\cite[Lem.~3.4.8 and 3.4.10]{CartGoulToin22} is needed for the second of these
methods to yield \req{growth}.  Since these methods have a better
$\epsilon$-order complexity, we now deduce that it is now
\emph{$o\big(\epsilon^{-3/2}\big)$ rather than
$\calO\big(\epsilon^{-3/2}\big)$ for finding first-order critical
points}, and \emph{$o\big(\epsilon^{-3}\big)$ rather than
$\calO\big(\epsilon^{-3}\big)$ to find second-order ones.}

\subsection{Adaptive regularization methods}

The case of adaptive regularization methods is quite similar to that of
trust-region algorithms. Again
\[
\{f_k\} =\{f(x_k)\},
\ms
\omega_k = \mbox{a $q$-th order criticality measure},
\]
$\kappa_d$ is an algorithm-specific constant and $\calS_k$ is the
index set of the ``successful iterations''. The bound \req{growth} is
now guaranteed by \cite[Lem~2.4.1 and 2.4.2]{CartGoulToin22}
with $\kappa_b=0$ and $\beta$ again depends on which type of critical
points are sought and the degree of derivatives used.  Because a
specific discussion of every case may quickly become cumbersome, we
only list, in Table~\ref{ar-t}, the algorithms of interest, pointers
to the relevant proofs, criticality order $q$ and associated refined
complexity bounds resulting from Theorem~\ref{result}.

\begin{table}[htbp]
  \begin{center}
    \begin{tabular}{|l|l|c|c|}
\hline
Algo. & Proof & Critic. order & Refined complexity \\
\hline
{\sf AR1}       & \cite[Th.~2.4.3]{CartGoulToin22} & 1rst & $o\big(\epsilon^{-2}\big)$ \\
{\sf AR2}       & \cite[Th.~3.3.4]{CartGoulToin22} & 1rst & $o\big(\epsilon^{-3/2}\big)$ \\
{\sf AR2}       & \cite[Th.~3.3.9]{CartGoulToin22} & 2nd  & $o\big(\epsilon^{-3}\big)$ \\
{\sf AR$p$}     & \cite[Th.~4.1.5]{CartGoulToin22} & 1rst, 2nd, 3rd & $o\big(\epsilon^{-(p+1)/(p+1-q)}\big)$ \\
{\sf AR$qp$}    & \cite[Th.~12.2.14]{CartGoulToin22} & 1rst, 2nd &$o\big(\epsilon^{-(p+1)/(p-q+1)}\big)$ \\
{\sf AR$qp$}    & \cite[Th.~12.2.14]{CartGoulToin22} & $q$-th, $q>2$ &$o\big(\epsilon^{-q(p+1)/p}\big)$ \\
{\sf AR$qp$IDA} & \cite[Th.~13.1.19]{CartGoulToin22} & 1rst, 2nd &$o\big(\epsilon^{-(p+1)/(p-q+1)}\big)$ \\
{\sf AR$qp$IDA} & \cite[Th.~13.1.19]{CartGoulToin22} & $q$-th, $q>2$ &$o\big(\epsilon^{-q(p+1)/p}\big)$ \\
{\sf AR$qp$EDA} & \cite[Th.~13.3.8]{CartGoulToin22} & 1rst, 2nd &$o\big(\epsilon^{-(p+1)/(p-q+1)}\big)$ \\
{\sf AR$qp$EDA} & \cite[Th.~13.3.8]{CartGoulToin22} & $q$-th, $q>2$ &$o\big(\epsilon^{-q(p+1)/p}\big)$ \\
{\sf AN2C}      & \cite[Th.~1]{GratJeraToin23a}& 1rst & $o\big(\epsilon^{-3/2}\big)$ \\
{\sf AN2C}      & \cite[Th.~2]{GratJeraToin23a}& 2nd & $o\big(\epsilon^{-3}\big)$ \\
{\sf AR$1p$GN}  & \cite[Th.~3.5]{GratToin23} & 1rst & $o\big(\epsilon^{-(p+1)/p}\big)$ \\
{\sf AR2GN}     & \cite[Th.~4.5]{GratToin23} & 1rst, 2nd & $o\big(\epsilon^{-(p+1)/(p+1-q)}\big)$ \\
\hline
    \end{tabular}
    \caption{\label{ar-t}Refined complexity bound for unconstrained adaptive
      regularization algorithms}
  \end{center}
\end{table}

The proofs for the AN2C algorithms \cite{GratJeraToin23a}, using an alternative
regularization of Newton's method, are more involved because $\calS_k$
is then the union of smaller sets, but again rely on ``telescoping
sums'' for subsets of iterations, \cite[Lem.~1 and 4]{GratJeraToin23a}
being used to ensure \req{growth} in this case. The {\sf AR$1p$GN} and
{\sf AR2GN} algorithms proposed in \cite{GratToin23} allows the use of
a general nonsmooth regularization, and \req{growth} is ensured by
\cite[Lem.~2.4 and 3.3]{GratToin23} in this case.

As it turns out, the proofs listed in Table~\ref{ar-t} are themselves
templates for the complexity proofs of variants of the adaptive
regularization that exploit problem structure. Again discussing every
case would be too cumbersome, but we refer the reader to
\cite{CartGoulToin22b} for a specialized algorithm for least Euclidean
distance optimization, to \cite[Th.~14.1.10]{CartGoulToin22} for a
variant designed for the minimization of possibly non-smooth composite
objectives, to \cite[Th.~13.1.19 and 13.3.8]{CartGoulToin22} (together
with \cite[Lem~13.1.9]{CartGoulToin22}) for noise-tolerant variants or
to \cite{CartGoulToin12a} for an algorithm exploiting
finite-differences approximations to derivatives, including the
derivative-free case.

\subsection{Direct-search methods}

Finally, direct-search methods for minimization may also be
considered.  In \cite[Th.~3.1]{Vice13} the telescoping sum argument is again
explicitly used to prove a worst-case complexity bound for this
important class of derivative-free methods. In this case, $\calS$ is
the set of successful iterations (as for trust-region and adaptive
regularization algorithms), $\{f_k\} =\{f(x_k)\}$, $\omega_k =
\|\nabla_x^1f(x_k)\|$, $\beta = 2$, $\kappa_d$ a constant involving the square
of the gradient's Lispchitz constant. The bound \req{growth} is
obtained by \cite[Lem.~3.2]{Vice13} for $k=k(\epsilon)$ (as needed to
derive \req{nearfinal}). The complexity bound in
$\calO\big(\epsilon^{-2}\big)$ of \cite[Cor.~3.1]{Vice13}
can then be refined to $o(\epsilon^{-2})$.

\subsection{Algorithms for constrained problems}

Because methods for unconstrained optimization do occur as crucial
ingredients of several algorithms for the constrained case, the
refined complexity bounds for the former may translate into refined
complexity bounds for the latter. The easiest situation is when
considering ``simple'' constraints, i.e.\ when the constraints define
a convex feasible set onto which projection is computationally affordable
(including, for example, the ubiquitous problem of minimizing a
function subject to simple bounds on the variables).  In this case, the
evaluation complexity bounds for unconstrained problems are often unmodified
(when considering their order as a function of $\epsilon$) compared
with their unconstrained counter-parts, and the techniques of proof to
establish them are directly derived from the unconstrained setting,
except for the use of criticality measures that are suitable for
constrained problems. See for instance \cite[Th.~6.2.3]{CartGoulToin22}, where the
$\calO(\cdot)$ bounds for first-, second- and third order critical
points may now be refined to $o(\cdot)$.

Finally, unconstrained or bound-constrained methods and the techniques
to prove their complexity are often instrumental in the analysis of
algorithms for more general nonlinear constraints (for instance for
``restoration'' or ``feasibility'' phases, where one minimizes the
violation of the nonlinear constraints, essentialy using algorithms
for unconstrained problems). Resulting bounds for the whole
constrained algorithm may then be refined along the lines described
above (see, for instance, \cite[Th.~7.2.2 and 7.2.6]{CartGoulToin22}
leading to \cite[Th.~7.2.7]{CartGoulToin22}).

\numsection{How close are the refined and standard bounds?}\label{ex-s}

Having discussed refined bounds for a significant selection of
algorithms, we now take a step back and investigate how much the
refined and standard bounds differ by looking at a particular example.
This example is univariate and built along the lines of
\cite[Th.~2.2.3]{CartGoulToin22} for steepest descent. Sequences of
iterates $\{x_k\}$, function values $\{f_k\}$, gradient values
$\{g_k\}$ and steps $\{s_k\}$ are first constructed to illustrate the
bound, and standard Hermite theory is then invoked to show the
existence of a suitable function interpolating these values.

\begin{figure}[htbp] 
  \centerline{
    \includegraphics[width=0.35\textwidth,height=4cm]{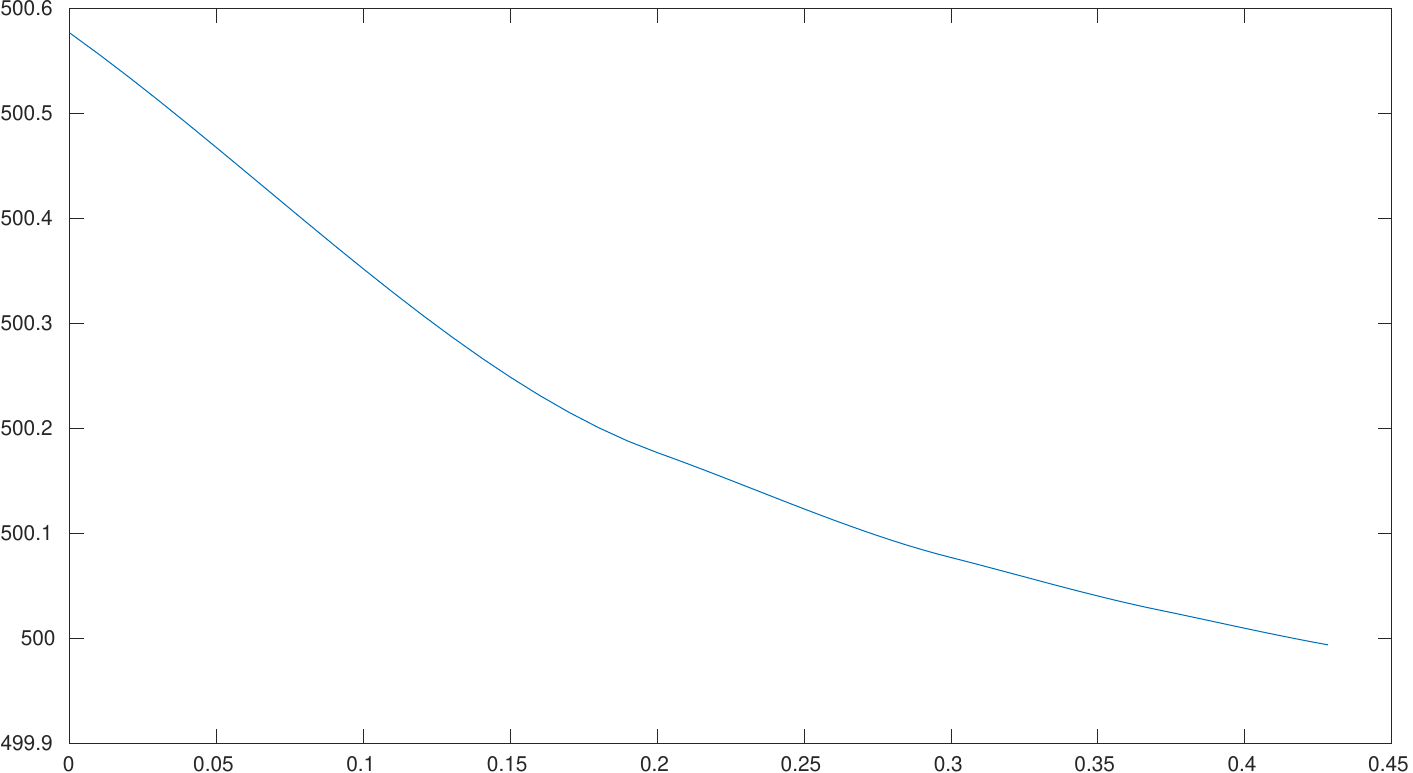}
    \hspace*{2mm}
    \includegraphics[width=0.35\textwidth,height=4cm]{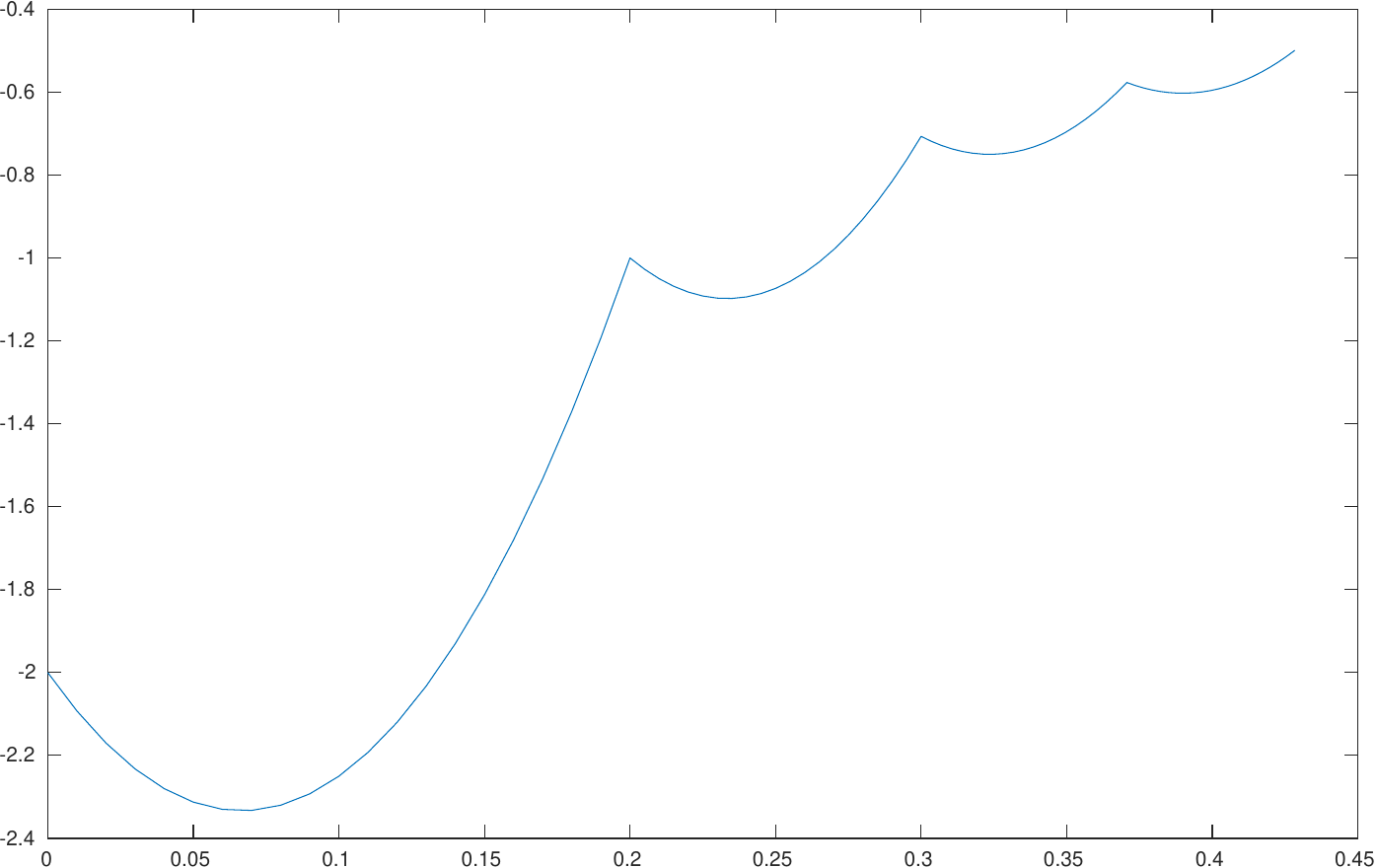}
    \hspace*{2mm}
    \includegraphics[width=0.35\textwidth,height=4cm]{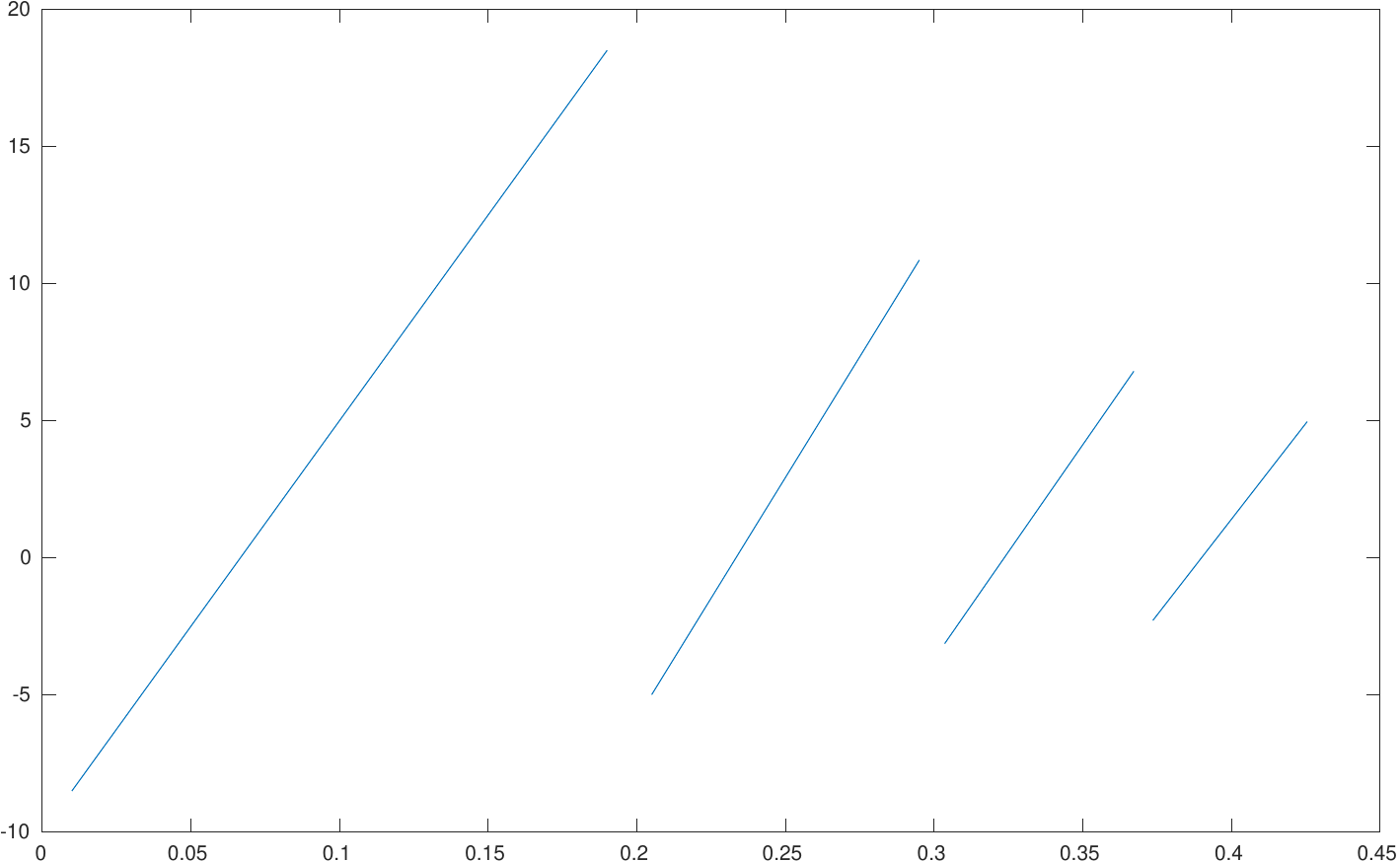}
  }
  \caption{\label{slow-f}The functions $f$ (left), $\nabla_x^1f$
    (middle) and $\nabla_x^2f$ (right) in $[x_0,x_{4}]$}
\end{figure}

Define, for $k\geq 0$ and some fixed constant $\delta>0$,
\beqn{gk-def}
g_0 = -2 \tim{ and } g_k = -\frac{1}{k^{\half+\delta}} \ms (k>0),
\eeqn
\beqn{fk-def}
f_0 = \zeta(1+2\delta)>1,
\ms
f_1 = f_0 - 4 \alpha
\tim{and}
f_{k+1} = f_k - \frac{\alpha}{k^{1+2\delta}} \ms (k>0)
\eeqn
and
\beqn{xk-def}
x_0 = 0,
\ms
x_1 = 2\alpha
\tim{and}
x_{k+1} = x_k + \frac{\alpha}{k^{\half+\delta}}  \ms (k>0).
\eeqn
for some $\alpha \in (0,1]$, where $\zeta(\cdot)$ is the Riemann function.
By definition of this function, we then have that $f_k\geq 0$ for all
$k\geq 0$, so the sequence $\{f_k\}$ is strictly decreasing and
bounded below, and hence convergent to some limit value $f_{\rm lim}\geq 0$.
As a consequence we have that
\beqn{limf}
\lim_{k \rightarrow \infty} (f_{\lfloor k/2 \rfloor} - f_k) = 0.
\eeqn
Moreover, let
\beqn{sk-def}
s_0 = 2\alpha
\tim{and}
s_k = \frac{\alpha}{k^{\half+\delta}}.
\eeqn
A simple calculation then shows that $k(\epsilon)$ as defined by
\req{keps-def} satisfies
\beqn{k-eps-def}
k(\epsilon)
= \lceil k_\epsilon \rceil
\tim{where}
k_\epsilon= \epsilon^{-\sfrac{1}{\half+\delta}}
= \epsilon^{-2} \epsilon^{\frac{4\delta}{1+2\delta}}
= o(\epsilon^{-2})
\eeqn
and thus
\beqn{tau-def}
k(\epsilon) \leq k_\epsilon +1= o(\epsilon^{-2}),
\eeqn
which is 
\req{bound2} (note that $k(\epsilon)$ tends to infinity
when $\epsilon$ tends to zero because of \req{gk-def}). But
\req{fk-def}, \req{k-eps-def} and \req{tau-def} together give that
\[
f_{\lfloor (k(\epsilon)-1)/2\rfloor} - f_{k(\epsilon)}
\geq\sum_{k=\lfloor (k(\epsilon)-1)/2\rfloor}^{k(\epsilon)-1}\frac{\alpha}{k^{1+2\delta}}
\geq \half |\calS_{k(\epsilon)-1}| \,\frac{\alpha}{k_\epsilon^{1+2\delta}}
= \half |\calS_{k(\epsilon)-1}| \,\frac{\alpha}{\epsilon^{-2}}
\]
and thus the stronger bound \req{bound} also holds. Now, because
\[
|f_{k+1}- f_k - g_ks_k| = 0 \leq s_k^2
\]
and
\[
|g_0-g_1| = |2-1| \leq \frac{1}{\alpha}s_0
\ms
|g_{k+1}-g_k| = \frac{1}{k^{\half+\delta}}-\frac{1}{(k+1)^{\half+\delta}}
\leq \frac{1}{\alpha}s_k,
\]
we may apply Hermite's theorem \cite[Th.~A.9.1]{CartGoulToin22} with
$\kappa_f=f_0$ on each interval $[x_k,x_{k+1}]$, defining a cubic
polynomial interpolating $f_k$, $f_{k+1}$, $g_k$ and
$g_{k+1}$. Combining these polynomials for successive intervals, we
obtain a cubic piecewise polynomial $f$ which is continuously differentiable
from $[0,+\infty)$ into $\Re$ and whose gradient is Lipschitz
continuous with a constant $L$ only depending on $\kappa_f>1$ and
$1/\alpha$. In addition, for all $k\geq 0$,
\[
f(x_k) = f_k
\tim{and}
\nabla_x f(x_k) = g_k
\]
and $f(x)$ is bounded below on $[0,+\infty)$. It is then easy to
extend this function on the left of the origin without altering
these properties by defining $f(x) = f_0-2x$ for $x<0$.
The left panel of Figure~\ref{slow-f} shows the (deceptively innocuous looking and barely nonconvex)
graph of $f$ in the interval $[x_0,x_4]$, where $\alpha
= 0.1$, $\delta = 0.001$ and $\epsilon=0.01$ (note that $\zeta(1.002) \approx 500.577$). The
middle panel shows the graph of its (continuous but non-monotone) gradient and the
right one its (discontinuous but bounded) Hessian.

As a consequence of the above argument, we may interpret the sequences
$\{x_k\}$, $\{f_k\}$ and $\{g_k\}$ as the result of a steepest descent
algorithm using an Armijo linesearch with initial stepsize $\alpha$,
applied to the univariate function $f$ and starting from $x_0$.  The
initial stepsize is always acceptable for the linesearch because
sufficient decrease ($\alpha g_k^2$) is achieved for the initial
stepsize at every iteration. In
view of \req{tau-def}, we have thus verified that, as expected, our
complexity bound in $o\big(\epsilon^{-2}\big)$ holds for the steepest
descent algorithm.

Now looking at \req{tau-def}, we also see that \textit{this bound can be
arbitrarily close to standard worst-case bound} in
$\calO(\epsilon^{-2})$ when $\delta$ is close to zero. A
similar conclusion also holds for the other cases discussed in
Section~\ref{algos-s}.

\numsection{Conclusion}

We have revisited the last step of the worst-case complexity proofs
for nonlinear optimization algorithms and obtained refined theoretical
bounds. We have then considered a few of the many cases where these
proofs can be refined, but the idea can clearly be applied more
widely. We have also shown that, although better, the refined bound
may be arbitrarily close to the standard one.

We note that our asymptotic results do not contradict the
non-asymptotic lower complexity bounds proved in
\cite{CartGoulToin10a} and \cite[Th.~2.2.3, 2.2.16 and
  12.2.17]{CartGoulToin22}.  Inded these latter bounds depend on
examples where a function is constructed such that convergence of the
relevant algorithm is exactly as slow as specified by the standard
$\calO(\cdot)$ bound. However, these functions explicitly depend on
$\epsilon$, which prevents taking the limit for $\epsilon$ tending to
zero, as we have done above. The situation is similar for the example
proposed in \cite{CarmDuchHindSidf20}, where the function on which
slow convergence occurs is defined in a space whose dimension depends
on $\epsilon$.

Of course, not all convergence proofs (and algorithms) are
concerned. Notable exceptions include complexity proofs for
measure-dominated problems (see \cite[Section~5.3]{CartGoulToin22})
because proofs in this context do not directly rely on the telescoping sum
argument.  The case of objective-function free (OFFO) algorithms, among which many
stochastic methods (see, for example,
\cite{GratRoyeViceZhan15,WardWuBott19,DefoBottBachUsun22,GratJeraToin24b}),
is less clear because the relevant complexity proofs typically involve telescoping sums along
with other potentially dominating terms. 

While the refined bounds are interesting, they remain of a fairly
generic nature, as we have verified in Section~\ref{algos-s}. It
remains an open question whether they can be refined further (maybe by
quantifying the numerator of the right-hand side of \req{bound}) for
specific methods.

{\footnotesize

  \section*{\footnotesize Acknowledgements}

  The authors wish to thank Oliver Hinder for an interesting
  discussion at ISMP 2024. The third author also gratefully acknowledges the
  friendly partial support of ANITI (Toulouse).


}

\end{document}